\documentclass[11pt]{article}
\usepackage{mathrsfs}
\usepackage{amscd}
\usepackage{amsmath,amsfonts,amssymb,amscd}
\usepackage{indentfirst,graphics,epsfig,psfrag}
\input{epsf}
\usepackage{ifpdf}
\usepackage{enumerate}
\usepackage{appendix}
\usepackage{enumerate}

\usepackage{lineno}

\makeatletter \@addtoreset{figure}{section} \makeatother
\makeatletter
\long\def\@makecaption#1#2{%
   \vskip 10\p@
   \setbox\@tempboxa\hbox{{#1}\ \ #2}%
   \ifdim \wd\@tempboxa >\hsize

       {#1}\ \ #2\par
   \else
       \hbox to\hsize{\hfil\box\@tempboxa\hfil}%
   \fi}
\makeatother

\newtheorem{thm}{Theorem}[section]

\newtheorem{lem}[thm]{Lemma}

\newtheorem{obs}[thm]{Observation}
\newtheorem{pro}[thm]{Proposition}

\newcommand{\qed}{{\hfill\rule{3pt}{7pt}}}

\def\qed{\hfill \rule{4pt}{7pt}}

\begin{document}
\title{\textbf{The generalized $3$-connectivity
of Lexicographic product graphs} \footnote{Supported by NSFC
No.11071130.}}
\author{
\small  Xueliang Li, Yaping Mao\\
\small Center for Combinatorics and LPMC-TJKLC\\
\small Nankai University, Tianjin 300071, P.R. China\\
\small E-mails: lxl@nankai.edu.cn; maoyaping@ymail.com
 }
\date{}
\maketitle
\begin{abstract}
The generalized $k$-connectivity $\kappa_k(G)$ of a graph $G$,
introduced by Chartrand et al., is a natural and nice generalization
of the concept of (vertex-)connectivity. In this paper, we prove
that for any two connected graphs $G$ and $H$, $\kappa_3(G\circ
H)\geq \kappa_3(G)|V(H)|$. We also give upper bounds for
$\kappa_3(G\Box H)$ and $\kappa_3(G\circ H)$. Moreover, all
the bounds are sharp.\\[2mm]
{\bf Keywords:} connectivity, Steiner tree, internally disjoint
trees, generalized connectivity, Lexicographic product.\\[2mm]
{\bf AMS subject classification 2010:} 05C05, 05C40, 05C76.
\end{abstract}

\section{Introduction}

All graphs considered in this paper are undirected, finite and
simple. We refer to the book \cite{bondy} for graph theoretical
notation and terminology not described here. The generalized
connectivity of a graph $G$, introduced by Chartrand et al. in
\cite{Chartrand1}, is a natural and nice generalization of the
concept of (vertex-)connectivity. For a graph $G=(V,E)$ and a set
$S\subseteq V(G)$ of at least two vertices, \emph{an $S$-Steiner
tree} or \emph{a Steiner tree connecting $S$} (or simply, \emph{an
$S$-tree}) is a such subgraph $T=(V',E')$ of $G$ that is a tree with
$S\subseteq V'$. Two $S$-trees $T$ and $T'$ are said to be
\emph{internally disjoint} if $E(T)\cap E(T')=\varnothing$ and
$V(T)\cap V(T')=S$. For $S\subseteq V(G)$ and $|S|\geq 2$, the
\emph{generalized local connectivity} $\kappa(S)$ is the maximum
number of internally disjoint $S$-trees in $G$. Note that when
$|S|=2$, an $S$-tree is just a path connecting the two vertices of
$S$. For an integer $k$ with $2\leq k\leq n$, the \emph{generalized
$k$-connectivity} $\kappa_k(G)$ of $G$ is defined as $\kappa_k(G)=
min\{\kappa(S) : S\subseteq V(G) \ and \ |S|=k\}$. Clearly, when
$|S|=2$, $\kappa_2(G)$ is nothing new but the connectivity
$\kappa(G)$ of $G$, that is, $\kappa_2(G)=\kappa(G)$, which is the
reason why one addresses $\kappa_k(G)$ as the generalized
connectivity of $G$. By convention, for a connected graph $G$ with
less than $k$ vertices, we set $\kappa_k(G)=1$. Set $\kappa_k(G)=0$
when $G$ is disconnected. Results on the generalized connectivity
can be found in \cite{Chartrand1, Chartrand2, LLMS, LLSun, LLL1,
LLL2, LL, LLZ, Okamoto}.

In addition to being a natural combinatorial measure, the
generalized connectivity can be motivated by its interesting
interpretation in practice. For example, suppose that $G$ represents
a network. If one considers to connect a pair of vertices of $G$,
then a path is used to connect them. However, if one wants to
connect a set $S$ of vertices of $G$ with $|S|\geq 3$, then a tree
has to be used to connect them. This kind of tree with minimum order
for connecting a set of vertices is usually called a Steiner tree,
and popularly used in the physical design of VLSI, see
\cite{Grotschel1, Grotschel2, Sherwani}. Usually, one wants to
consider how tough a network can be, for the connection of a set of
vertices. Then, the number of totally independent ways to connect
them is a measure for this purpose. The generalized $k$-connectivity
can serve for measuring the capability of a network $G$ to connect
any $k$ vertices in $G$.

Chartrand et al. {\upshape\cite{Chartrand2}} got the following
result for complete graphs.

\begin{thm}{\upshape\cite{Chartrand2}}\label{th1-1}
For every two integers $n$ and $k$ with $2\leq k\leq n$,
$\kappa_k(K_n)=n-\lceil k/2\rceil$.
\end{thm}

In \cite{LLZ}, Li et al. obtained the following results.

\begin{thm}{\upshape\cite{LLZ}}\label{th1-2}
For any connected graph $G$, $\kappa_3(G)\leq \kappa(G)$. Moreover,
the upper bound is sharp.
\end{thm}

\begin{thm} {\upshape\cite{LLZ}}\label{th1-3}
Let $G$ be a connected graph with at least three vertices. If $G$
has two adjacent vertices with minimum degree $\delta$, then
$\kappa_3(G)\leq \delta(G)-1$.
\end{thm}

\begin{thm}{\upshape\cite{LLZ}}\label{th1-4}
Let $G$ be a connected graph with $n$ vertices. For every two
integers $s$ and $r$ with $s\geq 0$ and $r\in \{0,1,2,3\}$, if
$\kappa(G)=4s+r$, then $\kappa_3(G)\geq 3s+\lceil\frac{r}{2}\rceil$.
Moreover, the lower bound is sharp.
\end{thm}

Sabidassi \cite{Sabidussi} derived the following perfect and
well-known theorem on the connectivity of Cartesian product graphs.

\begin{thm}{\upshape\cite{Sabidussi}}\label{th1-5}
Let $G$ and $H$ be two connected graphs. Then $\kappa(G\Box H)\geq
\kappa(G)+\kappa(H)$.
\end{thm}

Li et al. \cite{LLSun} studied the generalized $3$-connectivity of
Cartesian product graphs and got a lower bound of it. Their result
could be seen as an extension of Sabidussi' s Theorem.

\begin{thm}{\upshape\cite{LLSun}}\label{th1-6}
Let $G$ and $H$ be connected graphs such that
$\kappa_3(G)\geq\kappa_3(H)$. The following assertions hold:

$(i)$ If $\kappa(G)=\kappa_3(G)$, then $\kappa_3(G\Box H)\geq
\kappa_3(G)+\kappa_3(H)-1$. Moreover, the bound is sharp;

$(ii)$ If $\kappa(G)>\kappa_3(G)$, then $\kappa_3(G\Box H)\geq
\kappa_3(G)+\kappa_3(H)$. Moreover, the bound is sharp.
\end{thm}

From Theorem \ref{th1-5}, we know that $\kappa(G\Box H)\geq
\kappa(G)+\kappa(H)$. But we mention that it was incorrectly claimed
that $\kappa(G\Box H)=\kappa(G)+\kappa(H)$. In \cite{Spacapan},
\u{S}pacapan obtained the following result.

\begin{thm}{\upshape\cite{Spacapan}}\label{th1-7}
Let $G$ and $H$ be two nontrivial graphs. Then

$$
\kappa(G\Box
H)=\min\{\kappa(G)|V(H)|,\kappa(H)|V(G)|,\delta(G)+\delta(H)\}.
$$
\end{thm}

By the above result, we can derive a sharp upper bound of the
generalized $3$-connectivity for Cartesian product graphs.

\begin{thm}\label{th1-8}
Let $G$ and $H$ be two connected graphs. Then $\kappa_3(G\Box H)\leq
\min\{\lfloor
\frac{4}{3}\kappa_3(G)+r_1-\frac{4}{3}\lceil\frac{r_1}{2}\rceil\rfloor|V(H)|,
\lfloor\frac{4}{3}\kappa_3(H)+r_2-\frac{4}{3}\lceil\frac{r_2}{2}\rceil\rfloor
|V(G)|, \delta(G)+\delta(H)\}$, where $r_1\equiv \kappa(G) \ (mod \
4)$ and $r_2\equiv \kappa(H) \ (mod \ 4)$. Moreover, the bound is
sharp.
\end{thm}
\begin{pf}
From Theorem \ref{th1-4}, for a connected graph $G$, if
$\kappa(G)=4s+r_1$, then $\kappa_3(G)\geq
3s+\lceil\frac{r_1}{2}\rceil$, where $r_1\in \{0,1,2,3\}$. So
$\kappa_3(G)\geq 3\cdot
\frac{\kappa(G)-r_1}{4}+\lceil\frac{r_1}{2}\rceil=\frac{3}{4}
\kappa(G)-\frac{3}{4}r_1+\lceil\frac{r_1}{2}\rceil$, where
$r_1\equiv \kappa(G) \ (mod \ 4)$. Therefore, $\kappa(G)\leq
\frac{4}{3}\kappa_3(G)+r_1-\frac{4}{3}\lceil\frac{r_1}{2}\rceil$.
Similarly, for a connected graph $H$, $\kappa(H)\leq
\frac{4}{3}\kappa_3(H)+r_2-\frac{4}{3}\lceil\frac{r_2}{2}\rceil$,
where $r_2\equiv \kappa(H) \ (mod \ 4)$. From Theorem \ref{th1-2},
$\kappa_3(G\Box H)\leq \kappa(G\Box H)$. Furthermore, by Theorem
\ref{th1-5}, $\kappa_3(G\Box H)\leq \kappa(G\Box
H)=\min\{\kappa(G)|V(H)|,\kappa(H)|V(G)|,\delta(G)+\delta(H)\}\leq
\min\{\lfloor
\frac{4}{3}\kappa_3(G)+r_1-\frac{4}{3}\lceil\frac{r_1}{2}\rceil\rfloor|V(H)|,
\lfloor\frac{4}{3}\kappa_3(H)+r_2-\frac{4}{3}\lceil\frac{r_2}{2}\rceil\rfloor
|V(G)|, \delta(G)+\delta(H)\}$.\qed
\end{pf}

To show the sharpness of the above upper bound, we consider the
following example.

\noindent{\textbf{Example $1$.}}~Let $G=P_n \ (n\geq 3)$ and $H=P_m
\ (m\geq 3)$. Then $\kappa_3(G)=\kappa_3(H)=1$,
$\kappa(G)=\kappa(H)=1$ and hence $r_1=r_2=1$. From Theorem
\ref{th1-8}, $\kappa_3(G\Box H)\leq \min\{\lfloor
\frac{4}{3}\kappa_3(G)+r_1-\frac{4}{3}\lceil\frac{r_1}{2}\rceil\rfloor|V(H)|,
\lfloor\frac{4}{3}\kappa_3(H)+r_2-\frac{4}{3}\lceil\frac{r_2}{2}\rceil\rfloor
|V(G)|, \delta(G)+\delta(H)\}=\min\{n,m,2\}=2$. It can be checked
that for any $S\subseteq V(P_n\Box P_m)$ and $|S|=3$, $\kappa(S)\geq
2$, which implies $\kappa_3(G\Box H)\geq 2$. Thus, $P_n\circ P_m$ is
a sharp example for Theorem \ref{th1-8}.

Let us now turn our attention to another graph product. Recall that
the Lexicographic product of two graphs $G$ and $H$, written as
$G\circ H$, is defined as follows: $V(G\circ H)=V(G)\times V(H)$;
two distinct vertices $(u,v)$ and $(u',v')$ of $G\circ H$ are
adjacent if and only if either $(u,u')\in E(G)$ or $u=u'$ and
$(v,v')\in E(H)$. Let $G$ and $H$ be two connected graphs with
$V(G)=\{u_1,u_2,\ldots,u_{n}\}$ and $V(H)=\{v_1,v_2,\ldots,v_{m}\}$,
respectively. We use $G(u_j,v_i)$ to denote the subgraph of $G\circ
H$ induced by the set $\{(u_j,v_i)\,|\,1\leq j\leq n\}$. Similarly,
we use $H(u_j,v_i)$ to denote the subgraph of $G\circ H$ induced by
the set $\{(u_j,v_i)\,|\,1\leq i\leq m\}$. It is easy to see that
$G(u_{j_1},v_i)=G(u_{j_2},v_i)$ for different $u_{j_1}$ and
$u_{j_2}$ of $G$. Thus, we can replace $G(u_{j},v_i)$ by $G(v_i)$
for simplicity. Similarly, we can replace $H(u_{j},v_i)$ by
$H(u_j)$. For any $u,u'\in V(G)$ and $v,v'\in V(H)$, $(u,v),\
(u,v')\in V(H(u))$, $(u',v),\ (u',v')\in V(H(u'))$, $(u,v),\
(u',v)\in V(G(v))$, and $(u,v),\ (u',v)\in V(G(v))$. Note that
unlike the Cartesian Product, the Lexicographic product is a
non-commutative product. Thus $G\circ H$ needs not be isomorphic to
$H\circ G$.

In \cite{Yang}, Yang and Xu investigated the classical connectivity
of the Lexicographic product of two graphs.

\begin{thm}\label{th1-9}{\upshape\cite{Yang}}
Let $G$ and $H$ be two graphs. If $G$ is non-trivial, non-complete
and connected, then $\kappa(G\circ H)=\kappa(G)|V(H)|$.
\end{thm}

By the above result, we can derive a sharp upper bound of the
generalized $3$-connectivity of Lexicographic product graphs.

\begin{thm}\label{th1-10}
Let $G$ and $H$ be two connected graphs. If $G$ is non-trivial and
non-complete, then $\kappa_3(G\circ H)\leq \lfloor
\frac{4}{3}\kappa_3(G)+r-\frac{4}{3}\lceil\frac{r}{2}\rceil\rfloor
|V(H)|$, where $r\equiv \kappa(G) \ (mod~4)$. Moreover, the bound is
sharp.
\end{thm}
\begin{pf}
From Theorem \ref{th1-4}, for a connected graph $G$, if
$\kappa(G)=4s+r$, then $\kappa_3(G)\geq 3s+\lceil\frac{r}{2}\rceil$,
where $r\in \{0,1,2,3\}$. So $\kappa_3(G)\geq 3\cdot
\frac{\kappa(G)-r}{4}+\lceil\frac{r}{2}\rceil=\frac{3}{4}
\kappa(G)-\frac{3}{4}r+\lceil\frac{r}{2}\rceil$. Therefore,
$\kappa(G)\leq \lfloor
\frac{4}{3}\kappa_3(G)+r-\frac{4}{3}\lceil\frac{r}{2}\rceil\rfloor$.
From Theorem \ref{th1-2}, $\kappa_3(G\circ H)\leq \kappa(G\circ H)$.
Furthermore, by Theorem \ref{th1-9}, $\kappa_3(G\circ H)\leq
\kappa(G\circ H)=\kappa(G)|V(H)|\leq \lfloor
\frac{4}{3}\kappa_3(G)+r-\frac{4}{3}\lceil\frac{r}{2}\rceil\rfloor
|V(H)|$, where $r\equiv \kappa(G) \ (mod \ 4)$.\qed
\end{pf}

For the lower bound of $\kappa_3(G\circ H)$, we will prove the
following result next section.

\begin{thm}\label{th1-11}
Let $G$ and $H$ be two connected graphs. Then
$$
\kappa_3(G\circ H)\geq \kappa_3(G)|V(H)|.
$$
Moreover, the bound is sharp.
\end{thm}

To show the sharpness of the bounds of Theorems \ref{th1-10} and
\ref{th1-11}, we consider the following example.

\noindent{\textbf{Example $2$.}}~Let $G=P_n \ (n\geq 4)$ and
$H=P_3$. Then $\kappa_3(G)=\kappa_3(H)=1$. Since $\kappa(G)=1$, it
follows that $r=1$ and $\kappa_3(G\circ H)\leq \lfloor
\frac{4}{3}\kappa_3(G)+r-\frac{4}{3}\lceil\frac{r}{2}\rceil\rfloor|V(H)|=
3\lfloor \frac{4}{3}\kappa_3(G)-\frac{1}{3}\rfloor=3$ by Theorem
\ref{th1-10}. From Theorem \ref{th1-11}, we have $\kappa_3(G\circ
H)\geq \kappa_3(G)|V(H)|=3$. So $\kappa_3(G\circ H)=3$. Thus,
$P_n\circ P_3$ is a sharp example for both Theorem \ref{th1-10} and
Theorem \ref{th1-11}.

\section{Main results}

The following two lemmas are well known, which will be used later.

Given a vertex $x$ and a set $U$ of vertices, an \emph{$(x,U)$-fan}
is a set of paths from $x$ to $U$ such that any two of them share
only the vertex $x$.

\begin{lem}(Fan Lemma \cite{West}, p-170)\label{lem2-1}
A graph is $k$-connected if and only if it has at least $k+1$
vertices and, for every choice of $x$, $U$ with $|U|\geq k$, it has
an $(x,U)$-fan of size $k$.
\end{lem}

\begin{lem}(Expansion Lemma, \cite{West}, p-162)\label{lem2-2}
If $G$ is a $k$-connected graph, and $G'$ is obtained from $G$ by
adding a new vertex $y$ with at least $k$ neighbors in $G$, then
$G'$ is a $k$-connected.
\end{lem}

Let $G$ be a $k$-connected graph. Choose $U\subseteq V(G)$ with
$|U|=k$. Then the graph $G'$ is obtained from $G$ by adding a new
vertex $y$ and joining each vertex of $U$ and the vertex $y$. We
call this operation \emph{an expansion operation at $y$ and $U$}.
Denote the resulting graph $G'$ by $G'=G\vee \{y,U\}$.

We will prove Theorem \ref{th1-11} by three steps in the following
subsections.

\subsection{Lexicographic product of a path and a connected graph}

To start with, we show the following proposition.

\begin{pro}\label{pro2-3}
Let $H$ be a connected graph and $P_n$ a path with $n$ vertices.
Then $\kappa_3(P_n\circ H)\geq |V(H)|$. Moreover, the bound is
sharp.
\end{pro}

Let $H$ be a graph with $V(H)=\{v_1,v_2,\ldots,v_m\}$, and let
$V(P_n)=\{u_1,u_2,$ $\ldots,u_n\}$ such that $u_i$ and $u_j$ are
adjacent if and only if $|i-j|=1$. We need to show that for any
$S=\{x,y,z\}\subseteq V(P_n\circ H)$, there exist $m$ internally
disjoint $S$-trees. We proceed our proof by the following three
lemmas.

\begin{lem}\label{lem2-4}
If $x,y,z$ belongs to the same $V(H(u_i))$, $1\leq i\leq n$,
then there exist $m$ internally disjoint $S$-trees.
\end{lem}
\begin{pf}
Without loss of generality, we assume $x,y,z\in V(H(u_1))$. Then the
trees $T_i=x(u_2,v_i)\cup y(u_2,v_i)\cup z(u_2,v_i) \ (1\leq i\leq
m)$ are $m$ internally disjoint trees connecting $S$, as
desired.\qed
\end{pf}

\begin{lem}\label{lem2-5}
If only two vertices of $\{x,y,z\}$ belong to some copy $H(u_i) \
(1\leq i\leq n)$, then there exist $m$ internally disjoint
$S$-trees.
\end{lem}
\begin{pf}
We may assume $x,y\in V(H(u_1))$ and $z\in V(H(u_i)) \ (2\leq i\leq
n-1)$. In the following argument, we can see that this assumption
has no influence on the correctness of our proof. Consider the case
$i\geq 3$. Let $P'=u_2u_3\cdots u_n$. Clearly, $\kappa(P'\circ
H)\geq m$. From Lemma \ref{lem2-1}, there is a $(z,U)$-fan in
$P'\circ H$, where $U=V(H(u_2))=\{(u_2,v_j)|1\leq j\leq m\}$. Thus
there exist $m$ internally disjoint paths $P_1,P_2,\cdots,P_m$ such
that $P_j \ (1\leq j\leq m)$ is a path connecting $z$ and
$(u_2,v_j)$. Furthermore, the trees $T_j=x(u_2,v_j)\cup
y(u_2,v_j)\cup P_j \ (1\leq j\leq m)$ are $m$ internally disjoint
$S$-trees.

Now we assume $i=2$. We may assume $x,y\in V(H(u_1))$ and $z\in
V(H(u_2))$. Let $x',y'$ be the vertices corresponding to $x,y$ in
$H(u_2)$, $z'$ be the vertex corresponding to $z$ in $H(u_1)$.
Clearly, $H(u_1)$ is connected and so there is a path $P_1$
connecting $x$ and $y$ in $H(u_1)$.

If $z'\not\in \{x,y\}$, without loss of generality, let
$\{x,y,z'\}=\{(u_1,v_j)|1\leq j\leq 3\}$ and
$\{x',y',z\}=\{(u_2,v_j)|1\leq j\leq 3\}$, then the trees
$T_j=x(u_2,v_j)\cup y(u_2,v_j)\cup (u_1,v_j)(u_2,v_j)\cup z(u_1,v_j)
\ (4\leq j\leq m)$ and $T_1=xx'\cup x'y\cup yz$ and $T_2=xz\cup P_1$
and $T_3=xy'\cup y'z'\cup yy'\cup zz'$ are $m$ internally disjoint
$S$-trees; see Figure 2.1 $(a)$.

\begin{figure}[!hbpt]
\begin{center}
\includegraphics[scale=0.8]{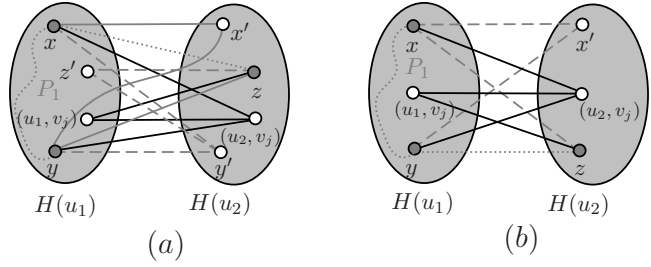}
\end{center}
\begin{center}
\caption{Graphs for Lemma \ref{lem2-5}.}
\end{center}\label{fig7}
\end{figure}

If $z'\in \{x,y\}$, without loss of generality, let $z'=y$,
$\{x,y\}=\{(u_1,v_1),(u_1,v_2)\}$ and
$\{x',z\}=\{(u_2,v_1),(u_2,v_2)\}$, then the trees
$T_j=x(u_2,v_j)\cup y(u_2,v_j)\cup (u_1,v_j)(u_2,v_j)\cup z(u_1,v_j)
\ (3\leq j\leq m)$ and $T_1=xz\cup xx'\cup yx'$ and $T_2=yz\cup P_1$
are $m$ internally disjoint $S$-trees; see Figure 2.1 $(b)$. The
proof is complete. \qed
\end{pf}

\begin{lem}\label{lem2-6}
If $x,y,z$ are contained in distinct $H(u_i)$s, then there exist $m$
internally disjoint $S$-trees.
\end{lem}
\begin{pf}
We have the following cases to consider.

{\flushleft\textbf{Case 1.}} $d_{P_n\circ H}(x,y)=d_{P_n\circ
H}(y,z)=1$.

We may assume that $x\in V(H(u_1))$, $y\in V(H(u_2))$, $z\in
V(H(u_3))$. In the following argument, we can see that this
assumption has no influence on the correctness of our proof. Let
$y',z'$ be the vertices corresponding to $y,z$ in $H(u_1)$, $x',z''$
be the vertices corresponding to $x,z$ in $H(u_2)$ and $x'',y''$ be
the vertices corresponding to $x,y$ in $H(u_3)$.

If $x,y',z'$ are distinct vertices in $H(u_1)$, without loss of
generality, let $\{x,y',z'\}=\{(u_1,v_j)|1\leq j\leq 3\}$ and
$\{x',y,z''\}=\{(u_2,v_j)|1\leq j\leq 3\}$ and
$\{x'',y'',z\}=\{(u_3,v_j)|1\leq j\leq 3\}$, then the trees
$T_j=x(u_2,v_j)\cup (u_1,v_j)(u_2,v_j)\cup y(u_1,v_j)\cup z(u_2,v_j)
\ (4\leq j\leq m)$ and $T_1=xx'\cup x'z'\cup x'z\cup z'y$ and
$T_2=xz''\cup zz''\cup y'z''\cup yy'$ and $T_3=xy\cup yz$ are $m$
internally disjoint $S$-trees; see Figure 2.2 $(a)$.

Assume that two of $x, y',z'$ are the same vertex in $H(u_1)$. If
$y'=z'$, without loss of generality, let
$\{x,y'\}=\{(u_1,v_1),(u_1,v_2)\}$ and $\{x',y\}=\{(u_2,v_1),\\
(u_2,v_2)\}$ and $\{x'',z\}=\{(u_3,v_1), (u_3,v_2)\}$, then the
trees $T_j=x(u_2,v_j)\cup (u_1,v_j)(u_2,v_j)\cup y(u_1,v_j)\cup
z(u_2,v_j) \ (3\leq j\leq m)$ and $T_1=xy\cup yz$ and $T_2=xx'\cup
x'x''\cup yx''\cup x'z$ are $m$ internally disjoint $S$-trees; see
Figure 2.2 $(b)$. The other cases ($x=y'$ or $x=z'$) can be proved
similarly.

\begin{figure}[!hbpt]
\begin{center}
\includegraphics[scale=0.8]{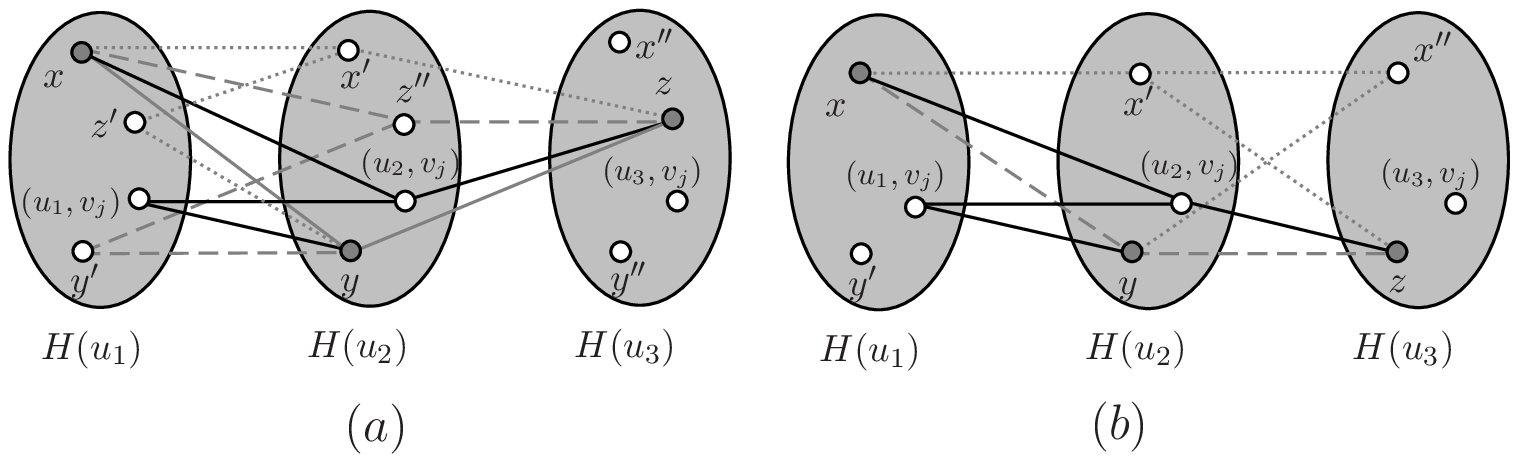}
\end{center}
\begin{center}
\caption{Graphs for Case $1$ of Lemma \ref{lem2-6}}.
\end{center}\label{fig7}
\end{figure}

Assume that $x,y',z'$ are the same vertex in $H(u_1)$. Without loss
of generality, let $x=(u_1,v_1)$, $y=(u_2,v_1)$ and $z=(u_3,v_1)$.
Then the trees $T_j=x(u_2,v_j)\cup (u_1,v_j)(u_2,v_j)\cup
y(u_1,v_j)\cup z(u_2,v_j) \ (2\leq j\leq m)$ and $T_1=xy\cup yz$ are
$m$ internally disjoint $S$-trees.

{\flushleft\textbf{Case 2.}} $d_{P_n\circ H}(x,y)=1$ and
$d_{P_n\circ H}(y,z)\geq 2$.

We may assume that $x\in V(H(u_1))$, $y\in V(H(u_2))$, $z\in
V(H(u_i)) \ (4\leq i\leq n)$. In the following argument, we can see
that this assumption has no influence on the correctness of our
proof. Let $y',z'$ be the vertices corresponding to $y,z$ in
$H(u_1)$, $x',z''$ be the vertices corresponding to $x,z$ in
$H(u_2)$ and $x'',y''$ be the vertices corresponding to $x,y$ in
$H(u_i)$. Let $P'=u_2u_3\cdots u_i$. Clearly, $\kappa(P'\circ H)\geq
m$. From Lemma \ref{lem2-1}, there is a $(z,U)$-fan in $P'\circ H$,
where $U=V(H(u_2))=\{(u_2,v_j)|1\leq j\leq m\}$. Thus there exist
$m$ pairwise internally disjoint paths $P_1,P_2,\cdots,P_m$ such
that each $P_j \ (1\leq j\leq m)$ is a path connecting $z$ and
$(u_2,v_j)$.

If $x,y',z'$ are distinct vertices in $H(u_1)$, without loss of
generality, let $\{x,y',z'\}=\{(u_1,v_j)|1\leq j\leq 3\}$ and
$\{x',y,z''\}=\{(u_2,v_j)|1\leq j\leq 3\}$, then the trees
$T_j=x(u_2,v_j)\cup (u_1,v_j)(u_2,v_j)\cup y(u_1,v_j)\cup P_j\
(4\leq j\leq m)$ and $T_1=xx'\cup x'y'\cup yy'\cup P_1$ and
$T_2=xy\cup P_2$ and $T_3=xz''\cup z'z''\cup z'y\cup P_3$ are $m$
internally disjoint $S$-trees; see Figure 2.3 $(a)$.

\begin{figure}[!hbpt]
\begin{center}
\includegraphics[scale=0.8]{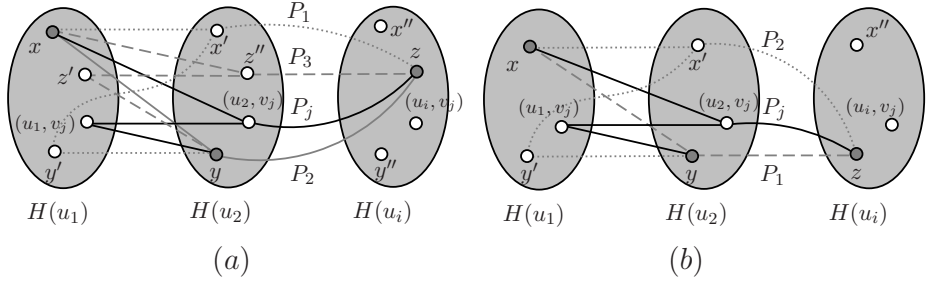}
\end{center}
\begin{center}
\caption{Graphs for Case $2$ of Lemma \ref{lem2-6}.}
\end{center}\label{fig7}
\end{figure}

Suppose that two of $x, y',z'$ are the same vertex in $H(u_1)$. If
$y'=z'$, without loss of generality, let
$\{x,y'\}=\{(u_1,v_1),(u_1,v_2)\}$ and $\{x',y\}=\{(u_2,v_1),\\
(u_2,v_2)\}$, then the trees $T_j=x(u_2,v_j)\cup
(u_1,v_j)(u_2,v_j)\cup y(u_1,v_j)\cup P_j\ (3\leq j\leq m)$ and
$T_1=xy\cup P_1$ and $T_2=xx'\cup x'y'\cup yy'\cup P_2$ are $m$
internally disjoint $S$-trees; see Figure 2.3 $(b)$. The other cases
($x=y'$ or $x=z'$) can be proved similarly.

Suppose that $x,y',z'$ are the same vertex in $H(u_1)$. Without loss
of generality, let $x=(u_1,v_1)$, $y=(u_2,v_1)$ and $z=(u_i,v_1)$.
Then the trees $T_j=x(u_2,v_j)\cup (u_1,v_j)(u_2,v_j)\cup
y(u_1,v_j)\cup P_j \ (2\leq j\leq m)$ and $T_1=xy\cup P_1$ are $m$
internally disjoint $S$-trees.

{\flushleft\textbf{Case 3.}} $d_{P_n\circ H}(x,y)\geq 2$ and
$d_{P_n\circ H}(y,z)\geq 2$.

We may assume that $x\in V(H(u_i))$, $y\in V(H(u_j))$, $z\in
V(H(u_k))$, where $i<j<k$, $|j-i|\geq 2$, $|k-j|\geq 2$, $1\leq
i\leq n-5$, $3\leq j\leq n-3$ and $5\leq k\leq n$. Let
$P'=u_i,u_{i+1},\cdots,u_{j-1}$ and
$P''=u_{j+1},u_{j+2},\cdots,u_{k}$. Then $P'$ and $P''$ are two
paths of order at least $2$. Since $\kappa(P'\circ H)\geq m$, from
Lemma \ref{lem2-2}, if we add the vertex $y$ to $P'\circ H$ and join
an edge from $y$ to each $(u_{j-1},v_r) \ (1\leq r\leq m)$, then
$\kappa((P'\circ H)\vee \{y,V(H(u_{j-1}))\})\geq m$. By the same
reason, $\kappa((P''\circ H)\vee \{y,V(H(u_{j+1}))\})\geq m$. From
Menger' s Theorem, there exist $m$ internally disjoint paths
connecting $x$ and $y$ in $(P'\circ H)\vee \{y,V(H(u_{j-1}))\}$, say
$P_1',P_2',\cdots,P_{m}'$. Also there exist $m$ internally disjoint
paths connecting $y$ and $z$ in $(P'\circ H)\vee
\{y,V(H(u_{j+1}))\}$, say $P_1'',P_2'',\cdots,P_{m}''$. Clearly,
each $P_i' \ (1\leq i\leq m)$ and each $P_j'' \ (1\leq j\leq m)$ are
also internally disjoint. Then the trees $T_i=P_i'\cup P_i'' \
(1\leq i\leq m)$ are $m$ internally disjoint $S$-trees. \qed
\end{pf}

From Lemmas \ref{lem2-4}, \ref{lem2-5} and \ref{lem2-6}, we conclude
that, for any $S\subseteq V(P_n\circ H)$, there exist $m$ internally
disjoint $S$-trees, namely, $\kappa(S)\geq m$. From the
arbitrariness of $S$, we have $\kappa_3(P_n\circ H)\geq m=|V(H)|$.
The proof of Proposition \ref{pro2-3} is complete.

\noindent{\textbf{Remark $1$.}}~As we have seen, for any
$S=\{x,y,z\}\subseteq V(P_n\circ H)$, there exist $m$ internally
disjoint $S$-trees in $P_n\circ H$. One can see that when $x,y,z$
belong to two copies $H(u_i)$ and $H(u_j)$ such that $u_iu_j\in
E(P_n) \ (1\leq i,j\leq n)$, we only use at most one path in
$H(u_i)$ or $H(u_j)$.

The sharpness of the bound in Proposition \ref{pro2-3} can be seen
from Example $2$.

\subsection{Lexicographic product of a tree and a connected graph}

In this subsection, we consider the generalized $3$-connectivity of
the Lexicographic product of a tree and a connected graph, which is
a preparation of the next subsection.

\begin{pro}\label{pro2-7}
Let $H$ be a connected graph and $T$ be a tree with $n$ vertices.
Then $\kappa_3(T\circ H)\geq |V(H)|$. Moreover, the bound is sharp.
\end{pro}
\begin{pf}
It suffices to show that for any $S=\{x,y,z\}\subseteq V(T\circ H)$,
there exist $m$ internally disjoint $S$-trees. Set
$V(T)=\{u_1,u_2,\ldots,u_n\}$, and $V(H)=\{v_1,v_2,\ldots,v_m\}$.

Let $x\in V(H(u_i)),y\in V(H(u_j)),z\in V(H(u_k))$ be three distinct
vertices. If there exists a path in $T$ containing $u_i$, $u_j$ and
$u_k$, then we are done from Proposition \ref{pro2-3}. If $i$, $j$
and $k$ are not distinct integers, such a path must exist. Thus,
suppose that $i$, $j$ and $k$ are distinct integers, and that there
exists no path containing $u_i$, $u_j$ and $u_k$. Then there exists
a subtree $T'$ in $T$ such that
$d_{T'}(u_i)=d_{T'}(u_j)=d_{T'}(u_k)=1$ and all the vertices of
$T'\setminus \{u_i,u_j,u_k\}$ have degree $2$ in $T'$ except for one
vertex, say $u_1$ with $d_{T'}(u_1)=3$. Clearly, there is a unique
path $P_1$ connecting $u_1$ and $u_i$, a unique path $P_2$
connecting $u_1$ and $u_j$, a unique path $P_3$ connecting $u_1$ and
$u_k$. Clearly, $P_1,P_2,P_3$ are pairwise internally disjoint. Let
$T''=T'\setminus \{u_i,u_j,u_k\}$.

If $d_{T'}(u_1,u_i)=d_{T'}(u_1,u_j)=d_{T'}(u_1,u_k)=1$, then the
trees $T_i=x(u_1,v_i)\cup y(u_1,v_i)\cup z(u_1,v_i) \ (1\leq i\leq
m)$ are $m$ internally disjoint $S$-trees.

Consider the case that $d_{T'}(u_1,u_i)\geq 2$,
$d_{T'}(u_1,u_j)=d_{T'}(u_1,u_k)=1$. It is clear that $T''\circ H$
is $m$-connected. Let $u_{i-1}$ be the vertex such that
$u_{i-1}u_i\in E(T')$ and $u_{i-1}$ is closer to $u_1$ than $u_{i}$
in $P_1$. From Lemma \ref{lem2-2}, $(T''\circ H)\vee
\{x,V(H(u_{i-1}))\}$ is $m$-connected and hence there exists an
$(x,U)$-fan in $(T''\circ H)\vee \{x,V(H(u_{i-1}))\}$, where
$U=V(H(u_1))$. So there exist $m$ internally disjoint paths
$P_{1,1},P_{1,2},\cdots,P_{1,m}$ connecting $x$ and
$(u_1,v_1),(u_1,v_2),\cdots,(u_1,v_m)$, respectively. Therefore, the
trees $T_i=P_{1,i}\cup y(u_1,v_i)\cup z(u_1,v_i) \ (1\leq i\leq m)$
are $m$ internally disjoint $S$-trees.

Consider the case that $d_{T'}(u_1,u_i)\geq 2$, $d_{T'}(u_1,u_j)\geq
2$ and $d_{T'}(u_1,u_k)=1$. Clearly, $T''\circ H$ is $m$-connected.
Let $u_{i-1}$ be the vertex such that $u_{i-1}u_i\in E(T')$ and
$u_{i-1}$ is closer to $u_1$ than $u_{i}$ in $P_1$. From Lemma
\ref{lem2-2}, $(T''\circ H)\vee \{x,V(H(u_{i-1}))\}$ is
$m$-connected and hence there exists an $(x,U)$-fan in $(T''\circ
H)\vee \{x,V(H(u_{i-1}))\}$. So there exist $m$ internally disjoint
paths $P_{1,1},P_{1,2},\cdots,P_{1,m}$ connecting $x$ and
$(u_1,v_1),(u_1,v_2),\cdots,(u_1,v_m)$, respectively (note that
$P_{1,1}$, $P_{1,2},\cdots,P_{1,m}$ belong to $P_1\circ H$).
Similarly, there exist $m$ internally disjoint paths
$P_{2,1},P_{2,2},\cdots,P_{2,m}$ connecting $y$ and $(u_2,v_1),
(u_2,v_2),\cdots,$ $(u_2,v_m)$, respectively  (note that $P_{2,1}$,
$P_{2,2},\cdots,P_{2,m}$ belong to $P_2\circ H$). Therefore, the
trees $T_i=P_{1,i}\cup P_{2,i}\cup z (u_1,v_i) \ (1\leq i\leq m)$
are $m$ internally disjoint $S$-trees.

Let us now consider the remaining case that $d_{T'}(u_1,u_i)\geq 2$,
$d_{T'}(u_1,u_j)\geq 2$ and $d_{T'}(u_1,u_k)\geq 2$. Similar to the
above method, there exist $m$ internally disjoint paths
$P_{1,1},P_{1,2},\cdots,P_{1,m}$ connecting $x$ and $(u_1,v_1),
(u_1,v_2),\cdots,(u_1,v_m)$, $m$ internally disjoint paths
$P_{2,1},P_{2,2},\cdots,P_{2,m}$ connecting $y$ and $(u_2,v_1),
(u_2,v_2)$, $\cdots,(u_2,v_m)$, and $m$ internally disjoint paths
$P_{3,1},P_{3,2},\cdots,P_{3,m}$ connecting $z$ and
$(u_3,v_1),(u_3,v_2),\cdots,(u_3,v_m)$, respectively. Therefore, the
trees $T_i=P_{1,i}\cup P_{2,i}\cup P_{3,i} \ (1\leq i\leq m)$ are
$m$ internally disjoint $S$-trees.

From the above argument, for any $S=\{x,y,z\}\subseteq V(G\circ H)$,
there exist $m$ internally disjoint $S$-trees, namely,
$\kappa(S)\geq m$. By the arbitrariness of $S$, we have
$\kappa_3(T\circ H)\geq m=|V(H)|$. The proof is complete.\qed
\end{pf}

To show the sharpness of Proposition \ref{pro2-7}, we consider the
following example.

\noindent{\textbf{Example $3$.}}~Let $H=P_3$, and $G=K_{1,n-1}$ be a
star with $n$ vertices, where $n\geq 4$. On one hand, from
Proposition \ref{pro2-7}, $\kappa_3(G\circ
H)=\kappa_3(K_{1,n-1}\circ P_3)\geq |V(P_3)|=3$. On the other hand,
let $V(G)=\{u_1,u_2,\cdots,u_{n}\}$ and $V(H)=\{v_1,v_2,v_3\}$. Let
$u_1$ be the center of the star $G=K_{1,n-1}$. Then
$u_2,u_3,\cdots,u_{n}$ are some leaves of this star. Choose
$S=\{(u_2,v_1),(u_3,v_1),(u_4,v_1)\}$. Clearly, there exist at most
three internally disjoint $S$-trees since
$V(H(u_1))=\{(u_1,v_1),(u_1,v_2),$ $(u_1,v_3)\}$. Therefore,
$\kappa_3(G\circ H)=\kappa_3(K_{1,n-1}\circ P_3)\leq 3$ and hence
$\kappa_3(P_n\circ P_3)=3$. Thus, $K_{1,n-1}\circ P_3$ is a sharp
example.

After the above preparations, we will prove our main result in the
following subsection.

\subsection{Lexicographic product of two general graphs}

Let $G$ and $H$ be two connected graphs. $G$ can be discomposed into
some trees $T_1,T_2,\cdots,T_k$. For each tree $T_i \ (1\leq i\leq
k)$, we can find some trees connecting $S\subseteq V(T_i\circ H)$
and $|S|=3$. In the end, we combine all the trees obtained by us,
and can obtain a lower bound of $\kappa_3(G\circ H)$ by the total
number of the trees connecting $S$.

\begin{figure}[!hbpt]
\begin{center}
\includegraphics[scale=0.7]{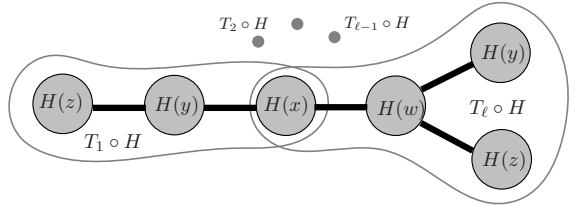}
\end{center}
\begin{center}
\caption{The structure of $(\bigcup_{i=1}^k T_i)\circ H$.}
\end{center}\label{fig7}
\end{figure}

\begin{obs}\label{obs2-8}
Let $G$ and $H$ be two connected graphs, $x,y,z$ be three distinct
vertices in $G$, and $T_1,T_2,\cdots,T_{\ell}$ be $\ell$ internally
disjoint $\{x,y,z\}$-trees in $G$. Then $(\bigcup_{i=1}^{\ell}
T_i)\circ H=\bigcup_{i=1}^{\ell} (T_i\circ H)$ has the structure as
shown in Figure $2.4$. Moreover, $(T_i\circ H)\cap (T_j\circ
H)=H(x)\cup H(y)\cup H(z)$ for $i\neq j$. In order to show the
structure of $(\bigcup_{i=1}^{\ell} T_i)\circ H$ clearly, we take
$\ell$ copies of $H(y)$, and $\ell$ copies of $H(z)$. Note that,
these $\ell$ copes of $H(y)$ (resp. $H(z)$) represent the same
graph.
\end{obs}

\noindent{\bf Example $4$.}~Let $H$ be the complete graph of order
$4$. The structure of $(T_1 \cup T_2)\circ H$ is shown in Figure
$2.5$.

\begin{figure}[!hbpt]
\begin{center}
\includegraphics[scale=0.7]{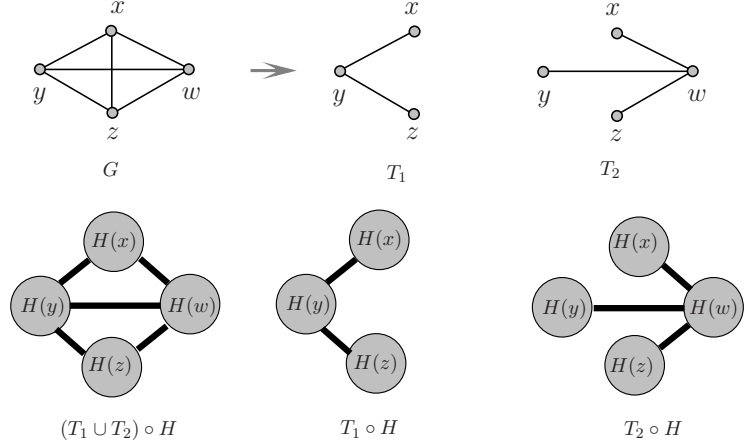}
\end{center}
\begin{center}
\caption{The structure of $(T_1\cup T_2)\circ H$.}
\end{center}\label{fig7}
\end{figure}

Now we are ready to prove Theorem $1.11$.

\textbf{Proof of Theorem $1.11$:} Without loss of generality, we set
$\kappa_3(G)=\ell$. It suffices to show that for any
$S=\{x,y,z\}\subseteq V(G\circ H)$, there exist $m\ell$ internally
disjoint $S$-trees. Assume that $V(G)=\{u_1,u_2,\ldots,u_n\}$ and
$V(H)=\{v_1,v_2,\ldots,v_m\}$.

Let $x\in V(H(u_i))$, $y\in V(H(u_j))$, $z\in V(H(u_k))$ be three
distinct vertices in $G\circ H$. Since $\kappa_3(G)=\ell$, there
exist $\ell$ internally disjoint $\{u_i,u_j,u_k\}$-trees $T_i \
(1\leq i \leq \ell)$ in $G$. By Observation \ref{obs2-8}, we know
that $(\bigcup_{i=1}^\ell T_i)\circ H=\bigcup_{i=1}^\ell (T_i\circ
H_i)$ and $(T_i\circ H_i)\cap (T_j\circ H_j)=H(x)\cup H(y)\cup H(z)$
for $i\neq j$.

Suppose that $x,y,z$ belong to the same $V(H(u_i))$, $1\leq i\leq
n$. Without loss of generality, let $x,y,z\in V(H(u_1))$. Since
$\delta(G)\geq \kappa_3(G)\geq \ell$, the vertex $u_1$ has $\ell$
neighbors in $G$, say $u_2,u_3,\cdots,u_{\ell+1}$. Then the trees
$T_{i,j}=x(u_i,v_j)\cup y(u_i,v_j)\cup z(u_i,v_j) \ (2\leq i\leq
\ell+1, 1\leq j\leq m)$ are $m\ell$ internally disjoint
$\{x,y,z\}$-trees in $G\circ H$.

Suppose that only two vertices of $\{x,y,z\}$ belong to some copy
$H(u_i) \ (1\leq i\leq n)$. Without loss of generality, let $x,y\in
H(u_1)$ and $z\in H(u_2)$. Since $\kappa(G)\geq \kappa_3(G)=\ell$,
there exist $\ell$ internally disjoint paths connecting $u_1$ and
$u_2$ in $G$, say $P_1,P_2,\cdots,P_{\ell}$. Clearly, there exists
at most one of $P_1,P_2,\cdots,P_{\ell}$, say $P_1$, such that
$P_1=u_1u_2$. From Remark $1$, there exist $m$ internally disjoint
$S$-trees in $P_1\circ H$, which occupies at most one path in
$H(u_1)$ or $H(u_2)$. For each $P_i \ (2\leq i\leq \ell)$, there
exist $m$ internally disjoint $S$-trees in $P_i\circ H$, which
occupies no edge in $H(u_j) \ (1\leq j\leq n)$. So the total number
of internally disjoint $S$-trees is $m\ell$, namely,
$\kappa_3(G\circ H)\geq m\ell$.

Assume that $x,y,z$ are contained in distinct $H(u_i)$s.  Without
loss of generality, let $x\in H(u_1)$, $y\in H(u_2)$ and $z\in
H(u_3)$. Since $\kappa_3(G)=\ell$, there exist $\ell$ internally
disjoint trees connecting $\{u_1,u_2,u_3\}$ in $G$, say
$T_1,T_2,\cdots,T_{\ell}$. For each tree $T_i \ (1\leq i\leq \ell)$,
there exist $m$ internally disjoint $S$-trees, which occupies no
edge in $H(u_j) \ (1\leq j\leq n)$. Thus, the total number of
internally disjoint $S$-trees is $m\ell$. So $\kappa_3(G\circ H)\geq
m\ell$.

From the above argument, we conclude that, for any $S\subseteq
V(G\circ H)$, $\kappa(S)\geq m\ell$, namely, $\kappa_3(G\circ H)\geq
m\ell=\kappa_3(G)|V(G)|$. The proof is complete. \qed

To show the sharpness of the bound in Theorem \ref{th1-11}, we can
also consider the following example.

\noindent{\textbf{Example $5$.}}~ Let $H=K_m$ be a complete graph
with $m$ vertices, and $G=P_n$ be a path with $n$ vertices, where
$n\geq 3$. On one hand, from Theorem \ref{th1-11}, $\kappa_3(G\circ
H)=\kappa_3(P_n\circ K_m)\geq m$. On the other hand,
$\kappa_3(G\circ H)\leq \kappa(G\circ H)=\kappa(P_n\circ K_m)=m$ by
Theorem \ref{th1-2}. So $\kappa_3(P_n\circ K_m)=m$. Thus, $P_n\circ
K_m$ is a sharp example.

\end{document}